\documentclass[reqno]{amsart}
\usepackage{amsfonts,amssymb,amsbsy,amsmath,amsthm,enumerate}
\newtheorem{theorem}{Theorem}[section]

\newtheorem{lemma}[theorem]{Lemma}

\newtheorem{definition}[theorem]{Definition}
\newtheorem{remark}[theorem]{Remark}

\newcommand{\bel}{\begin{equation} \label}
\newcommand{\ee}{\end{equation}}

\newcommand{\rd}{{\mathbb R}^{2}}
\newcommand{\re}{{\mathbb R}}

\newcommand{\la}{\langle}
\newcommand{\ra}{\rangle}

\DeclareMathOperator{\tr}{tr}

\newcommand{\w}{\omega}

\newcommand{\G}{\mathbb{G}}

\newcommand{\Q}{\mathcal{Q}}

\newcommand\R{\mathbb R}
\newcommand\N{\mathbb N}

\def\N{\mathbb{N}}
\def\P{\mathbb{P}}
\def\Z{\mathbb{Z}}

\def\E{\mathbb{E}}
\def\P{\mathbb{P}}

\def\I{\mathcal{I}}
\def\eps{\varepsilon}
\def\e{{\mathrm{e}}}

{

\begin{document}

\title[Splitting of the Landau levels and localization]
{Splitting of the Landau levels by magnetic perturbations and Anderson transition in 2D-random magnetic media}

\author[N. Dombrowski]{Nicolas Dombrowski}
\address[Dombrowski]{Universit\'e de Cergy-Pontoise,
CNRS UMR 8088, D\'epartement de Ma\-th\'e\-matiques, F-95000
Cer\-gy-Pontoise, France; {\em Present address}: Pontificia
Universidad Cat\'olica de Chile, Facultad de Matem\'aticas, Av.
Vicu\~na Mackenna 4860,  Santiago de Chile}
\email{nicolas.dombrowski@u-cergy.fr}

\author[F. Germinet]{Fran\c cois Germinet}
\address[Germinet]{Universit\'e de Cergy-Pontoise,
CNRS UMR 8088, IUF, D\'epartement de Math\'ematiques,
F-95000 Cergy-Pontoise, France}
\email{germinet@math.u-cergy.fr}

\author[G. Raikov]{Georgi Raikov}
\address[Raikov]{Pontificia Universidad Cat\'olica de Chile, Facultad de Matem\'aticas,
Av. Vicu\~na Mackenna 4860,  Santiago de Chile}
\email{graikov@mat.puc.cl}



\begin{abstract}
 {In this note we consider a Landau Hamiltonian perturbed by a random magnetic potential of Anderson type.
 For a given number of  bands, we prove the existence of both strongly localized states at the edges of the
 spectrum and dynamical delocalization near the center of the bands in the sense that wave packets travel at
 least at a given minimum speed. We provide explicit examples of magnetic perturbations that split the Landau
 levels into full intervals of spectrum.}
 \end{abstract}

\maketitle
\begin{center}
{\Large \em Dedicated to the memory of Pierre Duclos (1948 - 2010)} \\
\end{center}

\vspace{1cm}

 \today

\section{Introduction} \label{s1} \setcounter{equation}{0}

Over the past two decades  in the physics literature, a lot of
attention has been allocated to random magnetic fields in two
dimensions, see e.g. \cite{AHK,BSK,Fu,V} and the references therein.
The occurrence of localized states under the sole effect of a
random magnetic field has been recurrently predicted by the
theory, or computed at the band edges. It has then been a
challenging issue to provide evidence of the existence of extended
states in a 2D electron gas (2DEG) submitted to random magnetic
fields. It is commonly admitted that there are no such extended
states in 2DEG with random electric potentials in absence of a
constant perpendicular magnetic field (but no mathematical proofs
so far!). As far as random magnetic fields are concerned, the
issue is harder to settle for subtler effects seem to play a
role. For instance, while the computations of \cite{AHK} were
favorable to the occurrence of extended states, the ones from
\cite{BSK} were indicating their non existence.

In quantum Hall systems, namely a 2DEG submitted to a transverse
constant magnetic field, localized states are responsible for the
celebrated plateaux of the quantum Hall effect. In the case where
the Hall conductance is discontinuous, non trivial transport has
been proved to take place in \cite{GKS1} for electric disorder
(see also \cite{BES,GKS2,GKM}). In this note, we provide a similar
picture but with magnetic disorder. The random magnetic potential
is shown to create both strongly localized states at the edges of
the spectrum and dynamical delocalization near the center of the
band in the sense that wave packets travel at least at a given
minimum speed.

Mathematically, the proof of the occurrence of Anderson
localization due to random magnetic potentials only is not an easy
task, mainly because of the lack of monotonicity of the
eigenvalues as functions of the random variables. Very few
preliminary results are available: recently, Ghribi, Hislop, and
Klopp \cite{GHK} proved localization for  random magnetic
perturbations of a periodic magnetic potential in dimension $d \geq 2$
(see also \cite{KNNY} for a particular discrete model). They
exploit the Wegner estimate obtained in \cite{HK} together with
results of Ghribi \cite{Gh} in order to start a multiscale
analysis. In \cite{U}, Ueki extended \cite{HK} to prove
localization for some 2D-magnetic perturbation of the Landau
Hamiltonian at the very bottom of the spectrum (below the first
Landau level).

In this note we consider  2D-random magnetic perturbations of the
Landau Hamiltonian, and prove a transition between dynamical
localization and dynamical delocalization inside an arbitrary
number of bands. For our model, the phenomenon is thus similar to
that arising for random electric potentials \cite{GKS1,GKS2,GKM}.
The proof of localization exploits the Wegner estimate of Hislop
and Klopp \cite{HK}, revisited by Ghribi, Hislop, and Klopp
\cite{GHK}, together with a simple weak disorder argument to start
the multiscale analysis, provided some information on the location
of the spectrum that we address in a separate argument. Then
dynamical localization follows from \cite{GK1,GKgafa} together
with the full set of equivalent properties defining the region of
complete localization \cite{GKduke,GKjsp}. Delocalization is
proved along the lines of \cite{GKS1}; in particular the Hall
conductance is quantized, constant in the region of localization
and jumps by one as a Landau level is crossed. To our best
knowledge, this is the first 2D-random purely magnetic model for
which such a transition has been established mathematically.

If the theory of Anderson localization developed over the past
years for a continuum random Schr\"odinger applies, it remains to
prove that it does not lead to an empty result, namely that the
interval where states are shown to be localized does intersect the
spectrum. Getting detailed enough information about the location
of the almost spectrum of the random Hamiltonian is again trickier
with non monotonic perturbations of order 1 such as the magnetic
ones. In particular, in our setting, the issue reduces to the
proof that the Landau level does split as a periodic magnetic
perturbation is turned on and that the corresponding spectrum
contains an open interval. In \cite{DSS}, Dinaburg, Sinai, and
Soshnikhov considered small periodic electric perturbations of the
Landau Hamiltonian, and proved  that the low Landau levels split
into a set of positive Lebesgue measure. Gruber addressed the same
issue in \cite{Gr} but for magnetic perturbations. In this note,
we exhibit an explicit family of small periodic magnetic
perturbations for which the splitting gives rise to a full
interval of spectrum. This is achieved by a direct computation using
the translation invariance of our potential in one direction. Such
examples are then good enough to be randomized and used as random
magnetic fields.

This note is organized as follows. In Section~\ref{sectresult} we
introduce the model and state the main results. In
Section~\ref{sectsplit} we construct explicit examples for which
the low Landau levels split into  intervals as the magnetic
perturbation is turned on. In Section \ref{section4} we prove our
main result, Theorem \ref{thmresult}, while the Appendix, Section
\ref{section5}, contains some technical trace-class estimates used
in this proof.

\section{Main Results} \setcounter{equation}{0}
\label{sectresult}
    Let ${\bf A} = (A_1,A_2) \in L_{\rm loc}^2(\rd,
\rd)$ be a magnetic potential. Define the operator $H({\bf A})$ as
the self-adjoint operator generated in $L^2(\rd)$ by the closure
of the quadratic form
$$
\int_{\rd} |i\nabla u + {\bf A}u|^2 dx, \quad u \in C_0^{\infty}(\rd).
$$
The magnetic field generated by ${\bf A}$ is
$$
B(x) : = \frac{\partial A_2}{\partial x_1}(x)  - \frac{\partial
A_1}{\partial x_2}(x), \quad x = (x_1,x_2) \in \rd.
$$
 In the case of a constant magnetic field $B>0$ introduce the
 magnetic potential ${\bf A}_0 : = (0,Bx_1)$ which generates $B$. It is well-known
that the spectrum of $H({\bf A}_0)$  consists of the so-called
Landau levels $(2j-1)B$, $j \in {\mathbb N} : = \{1,2,\ldots\}$,
and each Landau level is an eigenvalue of $H({\bf A}_0)$ of
infinite multiplicity. The operator $H({\bf A}_0)$ is called the Landau Hamiltonian. From now on, the symbol $B$ stands for the magnetic strength of $H({\bf A}_0)$.

Let us introduce the random magnetic potential
    \bel{tt1}
{\bf A}_{\omega}(x)=  \sum_{\gamma\in\Z^2} \omega_\gamma {\bf
v}_\gamma(x),
    \ee
with ${\bf v}_\gamma(x)=(v_1(x -\gamma),v_2(x -\gamma))$, $\gamma
\in {\mathbb Z}^2$, $x \in \re^2$, $v_1,v_2$ being two given
$C^1(\R^2,\re)$ compactly supported functions, normalized so that
$\|\sum_{\gamma\in\Z^2} {\bf v}_\gamma\|_\infty=1$; the random
variables $(\omega_\gamma)_{\gamma\in\Z^2}$ are independent and
identically distributed, supported on $[-1,1]$, with common
density
$$
\rho_\eta(s)  = C_{\eta} \eta^{-1} \exp(-|s| \eta^{-1})
\chi_{[-1,1]}(s), \quad \eta>0,
$$
 and $C_\eta=(2(1-\e^{-\eta^{-1}}))^{-1}$ the normalization constant (note that $\frac12\le C_\eta \le 1$ for
$\eta\in]0,1]$). The support of $\rho_\eta$ is $[-1,1]$ for all
$\eta>0$, but as $\eta$ gets small, the disorder becomes weaker
in the sense that for most $\gamma$ the coupling $\omega_\gamma$
is  small. We may speak of a diluted random model (see \cite{GKS1,GKM} for a similar type of randomness). We denote by
$$
H_{B,\lambda,\omega,\eta} : = H({\bf A}_0 + \lambda {\bf
A}_{\omega,\eta}), \quad \lambda>0,
$$
the corresponding magnetic random operator, and will consider small values of the coupling constant $\lambda$.

For bounded Borelian functions $f$, the maps $\omega\to
f(H_{B,\lambda,\eta,\omega})$ are measurable. It follows from
standard ergodicity arguments that the spectrum is almost surely
deterministic as well as its pp, sc and ac components (see e.g.
\cite{KM,CFKS,St}). We denote by  $\Sigma_{B,\lambda}$ the almost sure spectrum of $H_{B,\lambda,\omega,\eta}$ (it does not depend on $\eta>0$ since by construction the support of $\rho_\eta$ is independent of $\eta>0$). It is easy to see that  $\Sigma_{B,\lambda}$ is contained in a union of intervals
$\I_j({B,\lambda})=[a_j(B,\lambda),b_j(B,\lambda)]\ni (2j-1)B$, $j
\in {\mathbb N}$. Moreover,  if\footnote{Here and in the sequel we write $a \lesssim b$ if there exists a constant $c$ such that $a \leq cb$.}
 ${\mathbb N} \ni J\lesssim
(B\lambda^2)^{-1}$, then
    \begin{equation}
\Sigma_{B,\lambda} \cap (-\infty, (2J-1)B+B] \subset
\bigcup_{j=1}^J \I_j({B,\lambda})\subset \bigcup_{j=1}^J
[(2j-1)B-C\lambda\sqrt{jB},
(2j-1)B+C\lambda\sqrt{jB}],\label{spectrum}
    \end{equation}
for some constant $C<\infty$ (see Lemma~\ref{lemgap} below). As a
consequence, for any integer $J \in {\mathbb N}$,  the first $J$
intervals $\I_j({B,\lambda})$, $j=1,\ldots,J$, are disjoint for
$\lambda$ small enough. More precisely,
  for any $B\in(0,\infty)$ there exists $\lambda_{\ast}$ such that for any
  $j\leq J$ and any $\lambda\in[0,\lambda_\ast)$ we have $\I_{j}({B,\lambda})\cap\I_{j+1}({B,\lambda})=\emptyset$,
  that is $b_{j}(B,\lambda) < a_{j+1}(B,\lambda)$. We denote by $\G_j(B,\lambda)
=(b_{j}(B,\lambda);a_{j+1}(B,\lambda))$ the $j$-th gap of the
spectrum. We say that the couple $(B,\lambda)$ respects the
\textit{ the disjoint band condition} if we have
\begin{equation}\label{DBC}
\G_j(B,\lambda)\neq\emptyset \text{ for any } j\leq J.
\end{equation}
It follows from \eqref{spectrum} that the disjoint band condition is satisfied if $\lambda\lesssim \sqrt{B/J}$.

\begin{definition}
The region of strong dynamical localization for $H_{B,\lambda,\omega,\eta} $ is denoted by $\Xi^{SDL}_{(B,\lambda,\eta)}\subset\R$, and is defined as the set of $E\in\R$ such that there exists an interval $I\ni E$ satisfying
\begin{equation}\label{dynloc}
\mathbb{E}\left\{ \sup_{t \in \R} \left\|  {\langle} x
{\rangle}^{\frac p 2} {\mathrm{e}^{-i tH_{B,\lambda,\eta,\omega} }}
\chi_I(H_{B,\lambda,\omega}) \tilde{\chi}_0 \right\|_2^2\right\}
   <\infty\
\end{equation}
for any $p>0$. Here $\|\cdot\|_2$ denotes the Hilbert-Schmidt norm, $\chi_I$ is
the characteristic function of $I$, and $\tilde{\chi}_0$ is the
characteristic function of the unit square centered at the origin.
\end{definition}

Strong dynamical localization is known to characterize the so called region of complete localization,
where many localization properties turn out to be equivalent  \cite{GKduke,GKjsp}. In particular, this region coincides with the set of energies
where the bootstrap multiscale analysis of \cite{GK1} applies.

\begin{theorem}\label{thmresult}
Fix $J \in {\mathbb N}$. Let $H_{B,\lambda, \omega,\eta}$ be the
Hamiltonian described above, satisfying the disjoint band condition \eqref{DBC}. Then there exists $\kappa_J>0$
(depending on $B$ and $J$) and $\Lambda=\Lambda(B,J)>0$, such that
for any $\lambda\in(0,\Lambda]$ and $\eta\in (0,c_{B,J}\lambda
|\log \lambda|^{-2}]$,  for all $j = 1,\cdots,J$, the Hamiltonian
$H_{B,\lambda,\eta,\omega}$ exhibits strong dynamical
localization, namely  for any interval $I$
satisfying
    \bel{oct16}
I\subset\Sigma_{B,\lambda}\cap [a_j(B,\lambda),
(2j-1)B-\kappa_J\lambda^2), \quad I\subset\Sigma_{B,\lambda}\cap
((2j-1)B+\kappa_J\lambda^2,b_j(B,\lambda)],
    \ee
    we have
    \begin{equation} \label{oct17}
I \subset \Xi^{SDL}_{(B,\lambda,\eta)} .
\end{equation}

Moreover for all $j=1,\cdots,J$, there exists a dynamical Anderson
transition $\tilde{E}_j(B,\lambda,\eta)\in
\Sigma_{B,\lambda}\cap[(2j-1)B-\kappa_J\lambda^2,
(2j-1)B+\kappa_J\lambda^2]$. More precisely, there exists at least
one energy $\tilde{E}_j(B,\lambda,\eta)\in \Sigma_{B,\lambda}\cap
[(2j-1)B-\kappa_J\lambda^2, (2j-1)B+\kappa_J\lambda^2]$, such that
for every non-negative $\mathcal{X}\in {C}^\infty_{0}
(\mathbb{R})$ with $\mathcal{X} \equiv 1$  on some open interval
containing $\tilde{ E}_j(B,\lambda) $, and for all  $p>0$,  we
have
\begin{equation}\label{momentgrowth}
\frac1{T} \int_0^{\infty} \mathbb{E}\left\{ \left\|  {\langle} x
{\rangle}^{\frac p 2} {\mathrm{e}^{-i tH_{B,\lambda,\eta,\omega} }}
\mathcal{X}(H_{B,\lambda,\eta,\omega}) \tilde{\chi}_0
\right\|_2^2\right\} {\mathrm{e}^{-\frac{t}{T}}} \,{\rm d}t
   \ge \
C_{p,\mathcal{X}} \, T^{\frac p4 - 6} \ ,
\end{equation}
for all  $T \ge 0$  with  $  C_{p,\mathcal{X}} > 0 $.
\end{theorem}

\begin{remark} \label{remoct2} A priori, it is not obvious that there
exist non empty intervals $I$ satisfying \eqref{oct16}. In the
proof of Theorem \ref{thmapplic} below we construct a family of
random potentials ${\bf A}_\omega$ and constant magnetic fields
$B$ for which such intervals exist, and provide an estimate of
their size.
\end{remark}

\begin{remark}
One may compare Theorem \ref{thmresult} to \cite[Corollary 2.4]{GKS1} where the only disorder parameter is the dilution coefficient $\eta$ (see \cite{GKM} as well). In  \cite[Corollary 2.4]{GKS1}, the random potential is an electric one, and, with notations of the present article, the authors consider the case $\lambda=1$ and $\eta$ small. Then localization is proved up to a distance $\mathcal{O}(\eta |\log \eta|)$ from the Landau levels, while here we get to a distance $\mathcal{O}(\lambda^2)=o(\eta^{2-\eps})$  for any $\eps>0$, if we take $\eta=c_{B,J}\lambda |\log \lambda|^{-2}$. Such a better bound is due to the combined effect of both parameters $\lambda$ and $\eta$.
\end{remark}


\section{Splitting of the Landau levels}
 \label{sectsplit}
 \setcounter{equation}{0}
Recall that $B$ denotes the constant scalar magnetic field generated by the magnetic potential ${\bf A}_0$ of the Landau Hamiltonian $H({\bf A}_0)$, and $\Sigma_{B,\lambda}$ denotes the almost sure
spectrum of the operator $H_{B,\lambda,\omega, \eta} = H({\bf A}_0 +
\lambda {\bf A}_{\omega})$.\\

\begin{theorem}\label{thmapplic}
Fix $J \in {\mathbb N}$. Then there exists random magnetic
potentials ${\bf A}_{\omega}$ of the form \eqref{tt1} with ${\bf
v} = (v_1,v_2)$ described explicitly in \eqref{oct15}, and
$\tilde{\kappa}_J > 0$ and $\tilde{\lambda}_J > 0$, such that for any $B$
in the set ${\mathcal M}_J = {\mathcal M}_J({\bf A}_{\omega})
\subseteq(0, \infty)$ described explicitly in \eqref{oct6}, we
have
    \bel{oct10}
 [(2j-1)B-\tilde{\kappa}_J\lambda,(2j-1)B+
\tilde{\kappa}_J\lambda]\subset \Sigma_{B,\lambda}, \quad
j=1,\cdots, J,
    \ee
    provided $\lambda \in (0,\tilde{\lambda}_J]$.
\end{theorem}
\begin{remark} \label{remoct3}
It follows from  \eqref{spectrum} and Theorem~\ref{thmapplic} that
the edges of the almost sure spectrum satisfy
$$
 \tilde{\kappa}_J(B) \lambda \le (a_j(B,\lambda) -
(2j-1)B), \quad (b_j(B,\lambda) - (2j-1)B) \le C\sqrt{jB} \lambda,
$$
for all
$j=1,\cdots,J$, all $\lambda\in (0,\tilde{\lambda}_J)$, and all $\eta>0$.
\end{remark}
\begin{remark} \label{remoct4}
The complement of the set ${\mathcal M}_J$ is always finite.
Moreover, Remark \ref{remmarch} below provides a simple sufficient condition that ${\mathcal M}_J = (0,\infty)$; 
informally, this is the generic condition that the periodic function $a$ appearing in \eqref{alabala} has sufficiently many non vanishing Fourier coefficients. 
\end{remark} 
\begin{remark} \label{remjuly}
In the construction of the potential ${\bf A}_{\omega}$ within the proof of Theorem \ref{thmapplic} we assume that the function $a$ appearing in \eqref{alabala} is given, and describe the set of admissible fields $B \in (0,\infty)$ for which Theorem \ref{thmapplic} holds true. Of course, we could  start with {\em an arbitrary} given $B$, and  construct afterwards(a family) of $1$-periodic $a \in C^1(\re; \re)$ for which Theorem \ref{thmapplic} is valid.
\end{remark}

{\em Proof of Theorem \ref{thmapplic}}: The main idea of the proof
 is to construct an appropriate periodic
magnetic potential ${\bf A}_{\rm per}$ such that for every
$\lambda>0$ any given Landau level splits into an interval of
positive length of the spectrum $\sigma(H({\bf A}_0 + \lambda {\bf
A}_{\rm per}))$
 of the operator $H({\bf A}_0 + \lambda {\bf A}_{\rm per})$. After
 that, using ideas of \cite{KM}, we show that
    \bel{oct11}
 \sigma(H({\bf A}_0 + \lambda {\bf A}_{\rm per})) \subset
 \Sigma_{B,\lambda},
    \ee
 which implies \eqref{oct10}.
 First, we construct ${\bf A}_{\rm per}$. Let $a \in C^{1}(\re;
\re)$ be a 1-periodic function whose derivative does not vanish identically. Evidently, $a$ is bounded on $\re$.
 Set
\begin{equation} \label{alabala}
{\bf A}_{\rm per}(x) = (0, a(x_1)), \quad x = (x_1,x_2) \in \rd.
\end{equation}
Let ${\mathcal F}$ be the partial Fourier transform with respect
to $x_2$, i.e.
$$
({\mathcal F}u)(x_1, k) : = (2\pi)^{-1/2} \int_{\re} e^{-ix_2 k}
u(x_1,x_2) dx_2.
$$
Then we have
$$
{\mathcal F} H({\bf A}_0 + \lambda {\bf A}_{\rm per}) {\mathcal
F}^* = \int_{\re}^{\oplus} \tilde{h}_{\lambda}(k) dk,
$$
where
    \bel{2}
 \tilde{h}_{\lambda}(k) : = - \frac{d^2}{dx_1^2} + (B x_1
+ \lambda a(x_1) - k)^2, \quad k \in \re,
    \ee
is the self-adjoint  operator in $L^2(\re)$,  essentially
self-adjoint on $C_0^{\infty}(\re)$. In \eqref{2} change the
variable $x_1 = t + k/B$, $t \in \re$. Then the operator
$\tilde{h}_{\lambda}(k)$ is unitarily equivalent to
$$
 h_{\lambda}(k) : = - \frac{d^2}{dt^2} + (B t
+ \lambda a(t + k/B))^2, \quad k \in \re.
    $$
    Note that we have
    \bel{3}
    h_{\lambda}(k) = h_{0} + v_{\lambda}(k),
    \ee
    where
    \bel{4}
    h_{0} : = - \frac{d^2}{dt^2} + B^2 t^2
    \ee
    and
    \bel{5}
    v_{\lambda}(t; k) : = 2 B t \lambda a(t + k/B) + \lambda^2 a(t +
    k/B)^2.
    \ee
    Since $B^2 t^2 + v_{\lambda}(t; k) \to \infty$ as $t \to \pm \infty$, the spectrum of
    the operator $h_{\lambda}(k)$ is discrete
    and simple.
    Let $\left\{E_j(k;\lambda)\right\}_{j=1}^{\infty}$ be the
    increasing sequence of its
    eigenvalues. Since the operators $h_{\lambda}(k)$ and $ \tilde{h}_{\lambda}(k)$ are unitarily equivalent, of course, $\left\{E_j(k;\lambda)\right\}_{j=1}^{\infty}$ is also the
    increasing sequence of the
    eigenvalues of $ \tilde{h}_{\lambda}(k)$. Applying an appropriate infinite-dimensional version of \cite[Theorem 5.16]{K} to the operator $ \tilde{h}_{\lambda}(k)$, we easily find
    that the functions $E_j$, $j \in
    {\mathbb N}$, are real analytic functions with respect to $k$ and $\lambda$ (see also the first footnote of \cite[p. 117]{K}). Moreover, $E_j(\cdot; \lambda)$, $j \in
    {\mathbb N}$,
    are $B$-periodic for any $\lambda \in \re$. It is well-known that
    \bel{15}
    \sigma(H({\bf A}_0 + \lambda {\bf A}_{\rm per})) = \bigcup_{j
    \in {\mathbb N}} \bigcup_{k \in [0,B)} \left\{E_j(k; \lambda)\right\}.
    \ee
    Further,
    \bel{12}
E_j(k;0) = (2j-1)B, \quad j \in {\mathbb N},
 \ee
i.e. the eigenvalues $E_j(k;0)$ coincide with the Landau levels
and are independent of $k \in \re$. Let   $\varphi_j$, $j \in
{\mathbb N}$, be the real eigenfunctions of the harmonic
oscillator $h_0$ which satisfy $h_0 \varphi_j = (2j-1)B \varphi_j$
and $\int_{\re} \varphi_j(t)^2 dt = 1$. We recall that
    \bel{6}
\varphi_j(t) = \varphi_j(t;B) = \frac{B^{1/4}}{\sqrt{(j-1)!
2^{j-1} \sqrt{\pi}}} {\mathcal H}_{j-1}(\sqrt{B} t) e^{-B t^2/2},
\quad j \in {\mathbb N}, \quad t \in \re,
    \ee
    where
   \bel{herm}
    {\mathcal H}_q(t) : = e^{t^2/2} \left(\frac{d}{dt} - t\right)^q
e^{-t^2/2}, \quad q \in {\mathbb Z}_+ : = \{0,1\ldots\},
    \ee
    are the Hermite polynomials.
 Fix $j \in {\mathbb N}$. Now, the so-called Feynman-Hellmann
    formula implies
    \bel{7}
    \frac{\partial E_j(k; 0)}{\partial \lambda} = 2B \int_{\re} a(t + k/B) t \varphi_j(t)^2 dt, \quad k
    \in \re.
    \ee
Assume that for some $k_{\pm}\in [0,B)$ we have
    \bel{7a}
    \frac{\partial E_j(k_-; 0)}{\partial \lambda} < 0, \quad \frac{\partial E_j(k_+; 0)}{\partial \lambda}>0.
    \ee

Taking into account relations \eqref{7a} and \eqref{12}, as well
as the continuity of $\frac{\partial E_j(k_{\pm};
\lambda)}{\partial \lambda}$ with respect to $\lambda$, we find
that there exist $\varkappa_j
> 0$ and $\lambda^*_j
> 0$ such that
    \bel{7b}
E_j(k_-;\lambda) - (2j-1)B < - \varkappa_j \lambda, \quad
E_j(k_+;\lambda) - (2j-1)B > \varkappa_j \lambda,
    \ee
    provided that $\lambda \in (0, \lambda^*_j)$. Combining
    \eqref{7} and \eqref{7a} with \eqref{7b} and \eqref{15}, we
    obtain the following
\begin{lemma} \label{newlemma1}
 Fix $B>0$ and $j \in {\mathbb N}$. Assume that for some
 $k_{\pm} \in [0,B)$ we have
    \bel{condsplit}
   \int_\R a(t+k_-/B) t \varphi_j(t)^2 dt < 0 ,  \quad \int_\R a(t+k_+/B) t \varphi_j(t)^2 dt > 0.
    \ee
    Then there exist
$\varkappa_j
> 0$ and $\lambda^*_j
> 0$ such that
    $$
[(2j - 1)B  - \varkappa_j \lambda,  (2j - 1)B + \varkappa_j
\lambda] \subset \sigma(H({\bf A}_0 + \lambda {\bf A}_{\rm per})),
    $$
    provided that $\lambda \in (0, \lambda^*_j)$.
\end{lemma}
Next we establish criteria which guarantee the existence of
$k_{\pm}$ for which inequalities \eqref{condsplit} hold true.\\
Expand $a$ into a Fourier series
\begin{equation}\label{decom}
a(t) = \sum_{l \in \mathbb Z} \alpha_l e^{i2\pi lt}, \quad t \in
\re,
\end{equation}
with Fourier coefficients
$$
\alpha_l = \overline{\alpha_{-l}} : = \int_0^1 a(t) e^{-i2\pi lt}
dt, \quad l \in {\mathbb Z}.
$$
Note that $a \in C^1(\re)$ implies $\{\alpha_l\}_{l \in {\mathbb
Z}} \in \ell^1(\mathbb Z)$. Moreover, by the assumption that the derivative of $a$ does not vanish identically,
there exists at least one non vanishing Fourier coefficient $\alpha_l$ with $l \in {\mathbb N}$. Further, we have
$$
F_j(k) = F_j(k;B) : = \int_\R a(t+k/B) t \varphi_j(t; B)^2 dt =
$$
    \bel{oct2}
    -2\sum_{l=1}^{\infty} |\alpha_l| I_j(2\pi l; B) \sin{\left(\frac{2\pi
    kl}{B} + \arg{\alpha_l}\right)}, \quad k \in \re, \quad j \in
    {\mathbb N},
    \ee
    where
    $$
    I_j(s; B) : = \int_{\re} \sin{(st)} t \varphi_j(t; B)^2 dt,
    \quad s \in \re.
    $$
    Evidently, $F_j$ is a $B$-periodic real analytic function of
    zero mean value. The existence of $k_{\pm} \in [0,B)$ for which inequalities \eqref{condsplit} hold
    true, is equivalent to the fact that $F_j$ does not vanish
    identically, which on its turn is equivalent to the existence
    of $l \in {\mathbb N}$ for which
    \bel{oct1}
    \alpha_l I_j(2\pi l; B) \neq 0.
    \ee
    \begin{remark} \label{remoct1} A condition which guarantees
    the splitting of the Landau levels into spectral bands of positive length, similar to
    \eqref{oct1}, was obtained in \cite[Chapter 4]{Be}, Note,
    however, that in \cite{Be} the Landau Hamiltonian perturbed by
    a periodic {\em electric} potential, was considered.
    \end{remark}
Let us now make condition \eqref{oct1} more explicit. Fix $j \in
{\mathbb N}$. Simple calculations yield
    \bel{oct3}
    I_j(s;B) = B^{-1/2} I_j(sB^{-1/2};1), \quad s \in \re,
    \ee
    and
    \bel{oct4}
    I_j(s;1) = \frac{-1}{(j-1)!2^{j-1} \sqrt{\pi}}\frac{d}{ds} \int_{\re} e^{ist} e^{-t^2} {\mathcal
    H}_{j-1}(t)^2dt, \quad s \in \re.
    \ee
    Applying \cite[Eq. (7.377)]{grry} (see also \cite[Lemma 2.2.2]{Be}), we get
    \bel{oct4a}
    \frac{1}{q!2^q \sqrt{\pi}} \int_{\re} e^{ist} e^{-t^2} {\mathcal
    H}_q(t)^2dt = {\mathcal L}_j(s^2/2) e^{-s^2/4}, \quad s \in
    \re, \quad q \in {\mathbb Z}_+
    \ee
    where
    $$
{\mathcal L}_q(\xi) : = \frac{1}{q!} e^{\xi} \frac{d^q}{d\xi^q}
(\xi^q e^{-\xi}), \quad \xi \in \re, \quad q \in {\mathbb Z}_+,
$$
are the Laguerre polynomials. Therefore,
    \bel{oct5}
    I_j(s;1) = {\mathcal P}_j(s) e^{-s^2/4}, \quad s \in \re,
    \ee
    where
    $$
    {\mathcal P}_j(s) : = \frac{s}{2} ({\mathcal L}_{j-1}(s^2/2) - 2 {\mathcal
    L}_{j-1}'(s^2/2)),
    $$
    and ${\mathcal L}_q'(\xi) = \frac{d{\mathcal L}_q(\xi)}{d\xi}$, $\xi \in
    \re$, $q \in {\mathbb Z}_+$.
    Obviously, ${\mathcal P}_j$, $j \in \N$, is an odd polynomial of degree $2j-1$. Hence, ${\mathcal P}_j$ has at most $j-1$ distinct positive real roots.
    \\
    Now  for $a\in C^1(\R;\R)$, recalling ~\eqref{decom}, set
    $$
    \mu_{j,l} : = \left\{s \in (0,\infty) \, | \, {\mathcal
    P}_j(2\pi ls) \neq 0\right\}, \quad l \in {\mathbb N}.
    $$
    Evidently, the complement in $(0,\infty)$ of $\mu_{j,l}$ contains no more than $j-1$ points. Put
    $$
    M_j : = \left\{\beta \in (0,\infty) \, | \, \beta^{-1/2} \in \bigcup_{l \in {\mathbb
    N}: \alpha_l \neq 0} \mu_{j,l}\right\}.
$$
Relations \eqref{oct3} -- \eqref{oct5} imply that the inclusion $B \in M_j$ holds if and only if at least one of the  coefficients
$$
\alpha_l I_j(2\pi l; B) = B^{-1/2} \alpha_l I_j(2\pi l B^{-1/2}; 1), \quad l \in \N,
$$
in \eqref{oct1} does not vanish, i.e. the inclusion $B \in M_j$ is equivalent to the existence of $k_\pm \in [0,B)$ for which \eqref{condsplit} is fulfilled. \\
Moreover, $B$ is in the complement of $M_j$, $j \in \N$, in $(0,\infty)$, if and only if ${\mathcal
    P}_j(2\pi l B^{-1/2}) = 0$ for each $l \in \N$ such that $\alpha_l \neq 0$. Thus, a simple
sufficient  (but not necessary) condition that $M_j = (0,\infty)$, is that the function $a$
has at least $j$ non vanishing Fourier coefficients $\alpha_l$
with $l \in {\mathbb N}$.\\
Thus, we obtain the following
\begin{lemma} \label{newlemma2}
 Fix $j \in {\mathbb N}$ and let $a\in C^1(\R;\R)$ be given. Then inequalities
    \eqref{condsplit}
   hold for some $k_{\pm} \in [0,B)$ if and only if $B \in M_j$.
\end{lemma}

Now we are in position to prove Theorem \ref{thmapplic}. Fix $J
\in {\mathbb N}$ and set
$$
\tilde{\lambda}_J : = \min_{j=1,\ldots,J} {\lambda}_j^*, \quad
\tilde{\kappa}_J = \min_{j=1,\ldots,J} \varkappa_j,
$$
and
    \bel{oct6}
    {\mathcal M}_J : = \bigcap_{j=1}^J M_j.
    \ee
    \begin{remark} \label{remmarch}
Evidently, the complement of the set ${\mathcal M}_J$ in
$(0,\infty)$ is always finite, and generically is empty. Similarly to $M_j$, a simple
sufficient  (but not necessary) condition that ${\mathcal M}_J = (0,\infty)$, is that the function $a$
has at least $J$ non vanishing Fourier coefficients $\alpha_l$
with $l \in {\mathbb N}$ (see \eqref{decom}).
\end{remark}

    Then Lemmas \ref{newlemma1} - \ref{newlemma2} imply that
    \bel{oct12}
    [(2j-1)B - \tilde{\kappa}_J \lambda, (2j-1)B + \tilde{\kappa}_J
    \lambda] \subset \sigma(H({\bf A}_0 + {\lambda
    A}_{\rm per})), \quad j=0,\ldots, J,
    \ee
    provided that $\lambda \in (0,\tilde{\lambda}_J)$. \\
Let $\zeta \in C_0^{\infty}(\re; \re)$ satisfy $0 \leq \zeta(t)
\leq 1$, $\sum_{m \in {\mathbb Z}} \zeta(t-m) = 1$, $t \in \re$.
Define ${\bf A}_\omega$ as in \eqref{tt1} with
    \bel{oct15}
v_1 = 0, \quad v_2(x) = a(x_1)\zeta(x_1)\zeta(x_2), \quad x =
(x_1, x_2) \in \rd.
    \ee
Evidently, ${\bf A}_\omega \in C^1(\rd; \rd)$, and
$$
\|{\bf A}_\omega\|_{L^\infty(\rd)} + \|\nabla {\bf
A}_\omega\|_{L^\infty(\rd)} < \infty,
$$
 for each
realization of the random variables $(\omega_\gamma)_{\gamma \in
{\mathbb Z}^2}$. Note that if $\tilde{\omega}$ is the periodic
realization of the random variables with $\tilde{\omega}_\gamma =
1$ for each $\gamma \in {\mathbb Z}^2$, then we have ${\bf
A}_{\tilde\omega} = {\bf A}_{\rm per}$, the magnetic potential
${\bf A}_{\rm per}$ being defined in \eqref{alabala}. Applying a
magnetic version of \cite[Theorem 4]{KM}, we find
    that \eqref{oct11} holds true. Finally, the combination of
    \eqref{oct11} and \eqref{oct12} yields \eqref{oct10}.

\section{ Proof of Theorem \ref{thmresult}}
\setcounter{equation}{0} \label{section4}

\subsection{First part: localization}
     To prove strong dynamical localization, we perform the bootstrap
multiscale analysis of \cite{GK1}. As it is well-known, multiscale
analysis in this context requires two main ingredients: a Wegner
estimate and an initial scale estimate. We shall play with small
enough $\lambda$'s and $\eta$'s to ensure these two ingredients. Since the
spectrum of $H_{B,\lambda,\eta,\omega}$ shrinks to the Landau
levels as $\lambda\to 0$, the constants appearing in the Wegner
estimate as well as in the multiscale analysis will grow as
$\lambda\to 0$. Since $\lambda$ will have to be chosen small
enough depending on those constants, precise versions of the
Wegner estimate and of the initial scale estimate are required.

Let ${\mathcal Q}_L \subset \rd$ be the square of side $L \in 6 {\mathbb N}$, centered
at the origin.
Let $\chi_L$ denote the characteristic function of ${\mathcal Q}_L$, and $\Gamma_L$
denote the characteristic function of the set ${\mathcal Q}_{L-1} \setminus {\mathcal Q}_{L-3}$.
Further, let $H_{B, \lambda, \eta,\omega, L}$ be the operator $(-i\nabla - {\bf A}_0 - \lambda {\bf A}_{\omega})^2$
with appropriate boundary conditions, self-adjoint in $L^2({\mathcal Q}_L)$.
For $z \in {\mathbb C} \setminus \sigma(H_{B, \lambda,\eta, \omega, L})$ set $R_{B, \lambda, \eta,\omega, L}(z) : = (H_{B, \lambda, \eta,\omega, L} - z)^{-1}$.

A Wegner estimate for random magnetic perturbation has been proved
in \cite[Thereom~6.1]{HK}. Because of the above considerations,
our analysis rather relies on the Wegner estimate obtained in
\cite{GHK}, and extended to the case of a random magnetic potential in \cite[Theorem~1.2]{HK}.
More precisely in our context \cite[Theorem~4.1]{GHK} reads:

\begin{theorem}\cite{GHK}\label{thmWegner}
Let $E\in ((2j-1)B,(2j+1)B)$, $j \in {\mathbb N}$, be given and set
$\delta={\rm dist}(E,\sigma(H({\bf A}_0)))$. Then there exists
$\lambda_0>0$ and, for any $q \in (0,1)$, a constant $Q_q<\infty$ such
that for any $\eps \in (0,\delta/2]$,  any $\lambda\le \lambda_0
\min\{1,\delta^{1/2}\}$, and  any $\eta>0$, we have
\begin{equation}\label{Wegner}
\P \left\{\mathrm{dist}
(E,\sigma(H_{B, \lambda,\eta, \omega, L}))\leq \eps
\right\}\le Q_W \eps^q L^q,
\end{equation}
with $Q_W=Q_q (\eta\delta)^{-1}$.
\end{theorem}

\begin{remark}
The factor $\eta^{-1}$ in \eqref{Wegner} comes from the derivative
of the probability distribution and \cite[Eq.~(3.16)]{HK} (see
also \cite[Theorem~1]{U}).
\end{remark}

For the initial scale estimate, we need the following version of
\cite[Theorem~2.4]{GKgafa}.

\begin{theorem}\cite[Theorem 2.4]{GKgafa}\label{thmfinvol}
Let $E\in ((2j-1)B,(2j+1)B)$, $j \in {\mathbb N}$, be given. Set
$\delta : ={\rm dist}(E,\sigma(H_{B}))$. Given a Wegner estimate of
the form \eqref{Wegner}, there exist $C_d,C_{d,q,j}<\infty$, so that
if for $L\ge C_d \delta^{-\frac3{16q} }$ we have
\begin{equation}\label{start}
\P\left(C_{d,q,j} B Q_W  L^{\frac{16}3}\| \Gamma_L
R_{B, \lambda, \omega, L}(E) \chi_{L/3}\| < 1\right) \ge 1 - 2.10^{-5},
\end{equation}
then $E\in\Xi^{SDL}_{(B,\lambda,\eta)}$.
\end{theorem}

In \eqref{start} we already took into account that the
constant  $\gamma_{\mathcal I}$ that appears in \cite[Assumption SLI]{GKgafa},  is
bounded by $C_d \sqrt{(2j+1)B}$. This can be seen from
\cite[Theorem~A.1]{GKduke}, since the magnetic perturbation is
relatively bounded with respect to $H_B$ with relative bound, say,
$\frac 12$.

We shall take advantage of the following lemma which is a
consequence of the resolvent identity (see \cite[Lemma~4.1]{DGR}).

\begin{lemma}\label{lemgap}
Let $H({\bf A}_0)$ be the Landau Hamiltonian with constant magnetic filed
$B$. Let ${\bf A}\in \mathcal{C}^1(\R^2)$ be such that $\|{\bf
A}\|_\infty\le K_1 \sqrt{B}$ and $\|\mathrm{div}\,{\bf A}\|_\infty
\le K_2 B$. Then there exists a constant $0<K_0<\infty$ such that
\begin{equation}
 \sigma(H({\bf A}_0+ {\bf A})) \cap [(2j-1)B - B,(2j-1)B +  B]
 \subset [(2j-1)B - d_j({\bf A},B),(2j-1)B +   d_j({\bf A},B)],
\end{equation}
for any $j$ so that where $ d_j({\bf A},B)< B$, where $d_j({\bf
A},B)=K_0\max(\|\mathrm{div}\,{\bf A}\|_\infty,\|{\bf
A}\|_\infty\sqrt{(j+1)B})$.

The same conclusions hold for the finite volume operators, with
the same constants, independent of the volume.
\end{lemma}

By the definition of the probability distribution,  given
$\eta\in]0,1]$, we have $\P(|\omega_0|\le \alpha) \ge 1 -
\exp(-\alpha\eta^{-1})$ (recall that the normalization constant of
the probability distribution satisfies $\frac12\le C_\eta \le 1$).
Now, for $B$ given and $\lambda\le \lambda_0$ (given by
Theorem~\ref{thmWegner}), we note  that the spectrum of
$H_{B, \lambda,\eta, \omega, L}$ satisfies
\begin{align}
&\mathbb{P}\left(\sigma(H_{B, \lambda,\eta, \omega, L})
\subset \bigcup_{j= 1}^J \left[B_j - C_j \lambda B^{1/2} \alpha, B_j + C_j \lambda B^{1/2} \alpha\right]\right) \\
& \qquad \ge \P\left(|\omega_j|\le \alpha, \,\forall  j\in\Lambda_{L}\right)\\
& \qquad\ge 1 -  \exp(-\alpha\eta^{-1}) L^2, \label{nospec}
\end{align}
with $0<\alpha<B^{1/2}$.

 Since we are working in spectral gaps,
we use the Combes-Thomas estimate of    \cite[Proposition
3.2]{BCH} (see also the proof of \cite[Theorem 3.5]{KK1} based on
\cite[Lemma 3.1]{BCH}), adapted to a finite volume as in
\cite[Section~3]{GKgafa}.

 Let
 $E \in \mathcal{I}_j(B)$ and assume that $|E-(2j-1)B|\ge2\delta$.
 We write $\delta=\kappa_J \lambda^2$, and choose $\kappa_J$ so that
 the condition $\lambda\le \lambda_0 \sqrt{\delta}$ in Theorem~\ref{thmWegner} is satisfied,
 namely $\kappa_J\ge \lambda_0^{-2}$. We further use \eqref{nospec} with $\alpha$ such that
 $C_J \lambda \sqrt{B} \alpha = \delta$, that is $\alpha = C_{B,J} \kappa_J \lambda$.

We pick $q\in (0,1)$ close to $1$, and $L$ such that $L>
\lambda^{-1} \ge C_J \lambda^{-3/(8q)}$ (hence, the assumption $L\ge
C_J \delta^{-3/(16q)}$ in Theorem~\ref{thmWegner} is fulfilled).

Then, using \eqref{nospec} and the Combes-Thomas estimate, we
conclude that condition \eqref{start} will be satisfied at the energy
$E$ if
\begin{gather}
 \alpha \eta^{-1}\ge C_3  \log L, \label{controlproba}\\
 C_{J,B} Q_q (\eta\delta)^{-1} {L}^{\frac {16} 3}
\e^{- C_4 \sqrt{\delta}{L}} < 1, \label{controlCT}
\end{gather}
where $C_3<\infty$ and $C_4 > 0$. Recalling that $\delta=\kappa_J
\lambda^2$, we choose $L/\log L \ge C_{B,J} \lambda^{-1} \log
(\lambda\eta)^{-1}$ so that \eqref{controlCT} holds, and
$\eta^{-1}\ge C_{B,J} \lambda^{-1}\log L$  so that
\eqref{controlproba} is satisfied. Since for $\eta$ small enough
$\eta^{-1}>>  \log \log \eta^{-1}$, these two conditions are
compatible, and $\eta\le c_{B,J}\lambda |\log \lambda|^{-2}$ is
sufficient.
As a consequence \eqref{start} holds for all $E$ so that $|E-(2j-1)B|\ge2\delta=2\kappa_J \lambda^2$, and strong dynamical localization follows.


\subsection{Second part: delocalization}

We finish the proof of  Theorem~\ref{thmresult} following the idea
of \cite{GKS1,GKS2} which consists in using the Hall conductance
in order to prove the existence of delocalization energies. In
\cite{GKS1,GKS2}, the authors considered a Landau Hamiltonian
perturbed by a random electric perturbation and proved various
properties of the Hall conductance, including the fact it is
integer valued in delocalization gaps. While both  \cite{GKS1} and
\cite{GKS2} can be generalized to the case of a random magnetic
perturbation, we focus on \cite{GKS1} for it provides a
simpler proof which does not require the more involved technology
of \cite{GKS2}.

Let $\Lambda$ be the characteristic function of the interval $[\frac{1}{2},\infty)$, and $\Lambda_j$ be  the multiplication operators given
by  $\Lambda_j(x)= \Lambda(x_j)$, $j=1,2$. For any orthogonal projection $P$ such that  $P[[P,\Lambda_1],[P,\Lambda_2]]$
is trace class, we set
\begin{equation}
\Theta(P):=\tr\{P[[P,\Lambda_1],[P,\Lambda_2]]\} = \tr [P\Lambda_1P,P\Lambda_2P] . \label{defTheta}
\end{equation}
Let $P$ be an orthogonal projection on $L^2(\re^2)$, and $\phi_x$
be a smooth characteristic function of the unit square centered at
$x \in \rd$.  Assume that we have
\begin{equation} \label{fin17}
\Vert \phi_x P \phi_y\Vert_2\leq K_p\la x \ra^\kappa \la y \ra^\kappa\e^{-\vert x-y\vert^\zeta}
\text{ for any } x, y \in\Z^2,
\end{equation}
with some $\zeta\in ]0,1]\ , \kappa>0$, and $K_P<\infty$. By
\cite[Lemma 3.1]{GKS1} we have
\begin{equation}
\vert \Theta \vert(P):=\Vert P[[P,\Lambda_1],[P,\Lambda_2]]
\Vert_1 \leq C_{\zeta, \kappa}K^2_P,
\end{equation}
where $\|\cdot\|_1$ is the trace-class norm. Then the so-called Hall
conductance is well-defined, and is given by
\begin{equation}
\sigma_{\mathrm{Hall}}(B,\lambda,\eta, E,\w)=-2\pi i \Theta(P_{B,\lambda,\eta, E,\w}).
\end{equation}
Applying the ergodic theorem  (see e.g. \cite{GKS1}), we obtain
\begin{equation}
\sigma_{\mathrm{Hall}}(B,\lambda,\eta, E):=\E \sigma_{\mathrm{Hall},\w}(B,\lambda,\eta,E,\w)=\sigma_{\mathrm{Hall}}(B,\lambda,\eta,E,\w) \text{ for } \P-\text{a.e } \w .
\end{equation}
We proceed as in \cite{GKS1} to get the existence of a delocalization energy near each Landau levels using a perturbative argument.

\begin{lemma}\label{int}
Assume that  $\lambda\lesssim \sqrt{B/J}$ so that the disjoint
band condition \eqref{DBC} holds. Then $\sigma_{\mathrm{Hall}}(B,\lambda,\eta,E)$
is constant in each connected component of
$\Xi_{(B,\lambda,\eta)}^{SDL}$. Moreover for any $j\le J$, we have
$\sigma_{\mathrm{Hall}}(B,\lambda,\eta, E)=j$ whenever
$E\in\Xi_{(B,\lambda,\eta)}^{SDL}\cap ] (2j-1)B, (2j+1)B[$.
\end{lemma}

\begin{proof}[Proof of Lemma \ref{int}]
That $\sigma_{\mathrm{Hall}}(B,\lambda,\eta,E)$ is constant in each connected
component of $\Xi_{(B,\lambda,\eta)}^{SDL}$ is a consequence of the
strong localization properties of the eigenfunctions that hold in
the region of strong dynamical localization. The  argument follows
from Lemma~3.1 and Lemma~3.2 of \cite{GKS1} which are general
results, independent of the particular form of the random
perturbation.

The proof of the second assertion is standard and consists in
starting with the zero disorder situation and a energy $E$ in the
middle of a given gap $\G_j(B,0)$, where  the Hall conductance
$\sigma_{\mathrm{Hall}}(B,\lambda,\eta,E)$ is known to be equal to $j$ (e.g.
\cite{ASS,BES}); then increase the disorder parameter $\lambda$
keeping $E$ in $\G_j(B,\lambda)$ and show that the conductance
remains constant; at last, use the fact that the Hall conductance
$\sigma_{\mathrm{Hall}}(B,\lambda,\eta,E)$ is constant when moving the energy $E$
inside a region of dynamical localization.

The first step, namely increasing the disorder, is a perturbative
argument which  is performed here for magnetic perturbations along the
lines of proof of \cite[Lemma~3.3]{GKS1}.

Pick $E=2jB$, the middle of the gap $\G_j(B,0)$. Since the gap
remains open for sufficiently small $\lambda \geq 0$, we can write
$P_\lambda=P_{B, \lambda,\eta,\w,E}$ as an appropriate Riesz
projection, apply the Combes-Thomas theory, and obtain the
estimate
$$
\Vert \phi_x P_\lambda \phi_y\Vert_2 \leq K_1 \e^{-K_1\vert x-y
\vert} \text { for all } x,y \in\Z^2 \text{ and } \lambda\in I,
$$
with some $K_1>0$ depending on $\eta$ (cf. \cite[Eq.
(3.16)]{GKS1}). In particular, \eqref{fin17} holds true.
     Suppose now that the perturbation ${\bf A}_\w$
has a compact support. By Lemma \ref{finl1}  the operator
$\mathcal{Q}_{\lambda,\lambda'}:=P_{\lambda}-P_{0} -
(P_{\lambda'}-P_{0})$  is trace-class for all $\lambda,
\lambda'\in I$. Using the second form of $\Theta(P_\lambda)$ in
\eqref{defTheta} and expanding the difference
$\Theta(P_\lambda)-\Theta(P_\lambda')$ in four terms with
$P_\lambda=P_{\lambda'} + \mathcal{Q}_{\lambda,\lambda'}$ as in
\cite[Eq. (3.35)]{GKS1} yields
$\Theta(P_\lambda)=\Theta(P_{\lambda'})$.

Next, we use an approximation argument considering
$\w^{L},\w^{>L}$ given by $\w_i^{L}=\w_i$ if $\vert i \vert \leq
L$ and $\w_i^{L}=0$ otherwise, and $\w_i^{>L}=\w_i-\w_i^{L}$ for
any $L>0$. With the obvious notations, we set
$\mathcal{Q}_{\lambda,>L}:=P_{\lambda}-P_{\lambda, L}$
 (cf.
\cite[Eq. (3.36)]{GKS1}). Using  an appropriate Combes-Thomas
estimate (see e.g. \cite[Lemma A.3]{CG}) and  a smooth partition of
unity $\{\phi_x\}_{x \in {\mathbb Z}^2}$, we find  that
    \bel{fin18}
\Vert \phi_x \Q_{\lambda,>L}\phi_y\Vert \leq C\e^{-C(\vert
x-y\vert+\max\{L-\vert x\vert,0\}+\max\{L-\vert y\vert,0\})}.
 \ee
Putting together \eqref{fin18} and \eqref{fin19}, we obtain the
estimate
\begin{align}
    \label{est}
    \Vert \phi_x \Q_{\lambda,>L}\phi_y\Vert_2
&\leq\Vert \phi_x \Q_{\lambda,>L}\phi_y\Vert^{\frac{1}{2}}
\Vert \phi_x \Q_{\lambda,>L}\phi_y\Vert^{\frac{1}{2}}_1\\
&\leq C'\e^{-C'(\vert x-y\vert+\max\{L-\vert
x\vert,0\}+\max\{L-\vert y\vert,0\})},
\end{align}
for all $x,y \in\Z^2$ and $L>0$.  Combining the fact that for any
orthogonal projections $P_\alpha, P_\beta, P_\gamma$ we have
$$
\Vert P_\alpha[[P_ \beta,\Lambda_1],[P_\gamma,\Lambda_2]] \Vert_1
\leq \sum\limits_{x,y,z\in\Z^2}\Vert \phi_x
[P_\beta,\Lambda_1]\phi_y \Vert_2 \Vert \phi_y
[P_\gamma,\Lambda_1]\phi_z\Vert_2 ,
$$
with \eqref{est}, and the dominated convergence theorem, we get
that $ \Theta(P_{\lambda})-\Theta(P_{\lambda,  L})\xrightarrow[L
\to \infty]{}0 $. This ends the proof of the lemma.
\end{proof}

We now finish the proof of the Theorem~\ref{thmresult}. Let us fix
the couple $(B,\lambda)$ so that  the disjoint band condition
\eqref{DBC} is valid. Pick $j\leq J$. By virtue of
Lemma~\ref{int}, it is not possible that
$\mathcal{I}_j(B,\lambda)\subset\Xi^{SDL}_{(B,\lambda,\eta)}$. As a
consequence, there exists at least one energy
$\tilde{E}_j(B,\lambda,\eta)\in\mathcal{I}_j(B,\lambda)$ such  that
$\tilde{E}_j(B,\lambda,\eta)\not\in\Xi^{SDL}_{(B,\lambda,\eta)}$. Next,
because of \eqref{oct16} and \eqref{oct17}, we have
$\tilde{E}_j(B,\lambda,\eta)\in[(2j-1)B-\kappa_J\lambda^2,
(2j-1)B+\kappa_J\lambda^2]$.

Finally, to get the dynamical lower bound \eqref{momentgrowth} near $\tilde{E}_j(B,\lambda,\eta)$ we apply \cite[Theorem 2.11]{GKduke}. We can indeed readily apply \cite{GKduke} to magnetic random perturbations: the bounds from \cite[Appendix A]{GKduke} are valid within our context since  $H_{B,\lambda,\omega,\eta} - H({\bf A}_0)$ is relatively $H({\bf A}_0)$--bounded with relative bound $<1$, and the proof of \cite[Theorem 2.11]{GKduke} itself  works as well for magnetic perturbations.


\section{Appendix: trace estimates}
\setcounter{equation}{0} \label{section5}
    Let ${\bf a} = (a_1,a_2) \in L_{\rm loc}^2(\rd; \rd)$. Introduce the self-adjoint
    operator $H_\lambda : = (-i\nabla - {\bf A}_0 - \lambda {\bf a})^2$,
    $\lambda \geq 0$,
    where, as earlier, the magnetic potential ${\bf A}_0$
    generates a constant magnetic field $B>0$. Denote by
    $P_{\lambda, E}$ the spectral projection of the
    operator $H_\lambda$ associated with the interval $(-\infty,
    E)$, $E \in \re$. We will say that ${\mathcal E} \in \re$ is
    in a spectral gap of the family $H_\lambda$, $\lambda \in
    [0,\lambda_0]$, with some $\lambda_0 > 0$, if there exist closed disjoint
    intervals ${\mathcal J}_-$ and ${\mathcal J}_+$ such that
    $$
    (-\infty, {\mathcal E}) \cap \bigcup_{\lambda \in [0,\lambda_0]}
    \sigma(H_\lambda) \subseteq {\mathcal J}_-, \quad
    ({\mathcal E}, \infty) \cap \bigcup_{\lambda \in [0,\lambda_0]}
    \sigma(H_\lambda) \subseteq {\mathcal J}_+.
    $$
    \begin{lemma} \label{finl1}
    Assume that ${\bf a} \in C^1(\rd; \rd)$ has a compact support. Let
    ${\mathcal E} \in \re$ is
    in a spectral gap of the family $H_\lambda$, $\lambda \in
    [0,\lambda_0]$ with some $\lambda_0 > 0$. Then the operator $P_{\lambda,
    {\mathcal E}} - P_{0,
    {\mathcal E}}$ is trace-class for all $\lambda \in
    [0,\lambda_0]$.
    \end{lemma}
    \begin{proof}
Evidently,
     there exists a
bounded contour $\Gamma$ such that ${\mathcal J}_-$ is contained
in its interior, ${\mathcal J}_+$ is contained in its exterior,
and there exists $s
> 0$ such that ${\rm dist}\,(\Gamma,
\sigma(H_\lambda)) > s$ for every $\lambda \in [0,\lambda_0]$. For
$z \in {\mathbb C}\setminus \sigma(H_\lambda)$ write $R_\lambda(z)
= (H_\lambda - z)^{-1}$. Then we have
    \bel{fin5a}
    P_{\lambda, {\mathcal E}} - P_{0, {\mathcal E}} =  \frac{1}{2\pi i} \int_\Gamma R_\lambda(z)
{\mathcal W} R_0(z) dz
    \ee
    where
    $$
    {\mathcal W} = {\mathcal W}_{\lambda} : = H_\lambda -  H_0 =
    2 \lambda {\bf a} \cdot (-i\nabla - {\bf A}_0) + i \lambda {\rm div}\,{\bf a} + \lambda^2
    |{\bf a}|^2.
    $$
Let $\zeta_0 \in C_0^{\infty}(\rd; \re)$ be a cut-off function,
equal to one on the support of ${\bf a}$. Then we have
  \bel{fin6}
    R_\lambda(z) {\mathcal W} R_0(z) =
     \zeta_0 R_\lambda(z) {\mathcal W} R_0(z) \zeta_0 +  R_\lambda(z) \zeta_1 S_\lambda
     R_\lambda(z) {\mathcal W} \zeta_0 R_0(z) + \zeta_0 R_\lambda(z)   {\mathcal W} R_0(z)  S_0 \zeta_1 R_0(z)
    \ee
    where
    $$
    S_\lambda = [H_\lambda(z),  \zeta_0] : =  2i \nabla \zeta_0 \cdot (-i\nabla - {\bf A}_0 - \lambda {\bf a}) - \Delta \zeta_0, \quad \lambda \geq 0,
    $$
    and $\zeta_1 \in C_0^{\infty}(\rd; \re)$ is a cut-off function, equal to one on the support of $\zeta_0$. Obviously,
  \bel{fin7}
    \left\|\int_\Gamma R_\lambda(z) \zeta_1 S_\lambda R_\lambda(z) {\mathcal W} \zeta_0 R_0(z) dz \right\|_1 \leq
   |\Gamma|  \sup_{z \in \Gamma} \left(\|R_\lambda(z) \zeta_1\|_2 \|S_\lambda R_\lambda(z) {\mathcal W}\| \|\zeta_0 R_0(z)\|_2\right),
     \ee
     \bel{fin8}
    \left\|\int_\Gamma \zeta_0 R_\lambda(z)   {\mathcal W} R_0(z)  S_0 \zeta_1 R_0(z) dz \right\|_1 \leq
     |\Gamma| \sup_{z \in \Gamma} \left(\|\zeta_0 R_\lambda(z) \|_2 \|{\mathcal W} R_0(z)  S_0 \| \|\zeta_1 R_0(z)\|_2\right),
     \ee
     where $\| \cdot \|_1$ denotes the trace-class norm, and $|\Gamma|$ is the length of  $\Gamma$.
     Applying the Hilbert-Schmidt diamagnetic inequality (see e.g. \cite[Theorem 2.13]{Si}), to the operators
     $\zeta_j R_\lambda(-1)$, we  find that
     $$
     \|\zeta_j R_\lambda(z)\|_2^2 = \| R_\lambda(z) \zeta_j\|_2^2 \leq
     \frac{c_0^2}{(2\pi)^2} \|\zeta_j\|^2_{L^2(\rd)} \int_{\rd} \frac{d\xi}{(|\xi|^2 + 1)^2},
     \quad
     j=0,1, \quad z \in \Gamma, \quad \lambda \in [0,\lambda_0],
     $$
     with
     $$
     c_0 : =  \sup_{z \in \Gamma} \sup_{\lambda \in [0,\lambda_0]}
     \sup_{E \in \sigma(H_\lambda)} \frac{E+1}{|E-z|}.
     $$
    Similarly,
      $$
      \sup_{z \in \Gamma} \|S_\lambda R_\lambda(z) {\mathcal W}\| <
      \infty,\quad
      \sup_{z \in \Gamma} \|{\mathcal W} R_0(z)  S_0
      \| < \infty.
      $$
Further, set $B_j : = B(2j-1)$, $j \in {\mathbb N}$, and
     write $R_0(z) = \sum_{j \in {\mathbb N}}(B_j - z)^{-1} \Pi_j$ where $\Pi_j$
     is the orthogonal projection onto ${\rm Ker}\,(H_0 - B_j)$. Put
     $$
     R_\lambda^{-}(z) : = \int_{(-\infty, {\mathcal E})} (E-z)^{-1} d_E P_{\lambda, E},
     \quad z \in {\mathbb C}\setminus (-\infty, {\mathcal E}),
     $$
     $$
     R_\lambda^{+}(z) : = \int_{({\mathcal E},\infty)} (E-z)^{-1} d_E P_{\lambda, E},
     \quad z \in {\mathbb C}\setminus ({\mathcal E},\infty).
     $$
     By the Cauchy theorem,
     \bel{fin1}
     \frac{1}{2\pi i} \int_\Gamma \zeta_0 R_\lambda(z) {\mathcal W} R_0(z) \zeta_0 dz  =
     \sum_{j \in {\mathbb N}\, : \, B_j\in {\mathcal J}_+} \zeta_0
     R_{\lambda}^-(B_j) {\mathcal W} \zeta_0 \Pi_j \zeta_0 -  \sum_{j \in {\mathbb N}\, : \, B_j\in {\mathcal J}_-}
     \zeta_0 R_{\lambda}^+(B_j) {\mathcal W} \zeta_0 \Pi_j \zeta_0.
     \ee
     Let us estimate the trace-class norm of the first (infinite) sum at the r.h.s. of \eqref{fin1}.
     For each $j \in {\mathbb N}$ such that $B_j\in {\mathcal J}_+$ we have
     $$
     \zeta_0 R_{\lambda}^-(B_j) {\mathcal W} \zeta_0 \Pi_j \zeta_0 =
     $$
     $$
     \zeta_0 R_{\lambda}^-(B_j) (H_\lambda + 1)^2 R_\lambda(-1)
     \zeta_0 R_\lambda(-1) {\mathcal W} \zeta_0 \Pi_j \zeta_0 +
     \zeta_0 R_{\lambda}^-(B_j) (H_\lambda + 1)^2 R_\lambda(-1)
     \zeta_1 S_\lambda R_\lambda(-1)^2 {\mathcal W} \zeta_0 \Pi_j \zeta_0.
     $$
     Therefore,
     $$
     \left\|\sum_{j \in {\mathbb N}\, : \, B_j\in {\mathcal J}_+} \zeta_0 R_{\lambda}^-(B_j) {\mathcal W}
     \zeta_0 \Pi_j \zeta_0\right\|_1 \leq
     $$
    \bel{fin2}
     \left(\|R_\lambda(-1) \zeta_0 \|_2 \|R_\lambda(-1) {\mathcal W}\| + \|R_\lambda(-1)
     \zeta_1 \|_2 \|S_\lambda R_\lambda(-1)^2 {\mathcal W}\|\right)
     \,
     \sum_{j \in {\mathbb N}\, : \, B_j\in {\mathcal J}_+} \|\zeta_0 R_{\lambda}^-(B_j)
     (H_\lambda + 1)^2\| \|\zeta_0 \Pi_j \zeta_0\|_2.
     \ee
     By the spectral theorem,
     \bel{fin3}
     \|\zeta_0 R_{\lambda}^-(B_j) (H_\lambda + 1)^2\| \leq \|\zeta_0\|_{L^{\infty}(\rd)} \sup_{E \in {\mathcal J}_-} \frac{(E+1)^2}{|E-B_j|} \leq c_1 j^{-1}
     \ee
     where $c_1$ is independent of $j$. Next, \cite[Lemma 3.1]{KP} implies
     \bel{fin4}
     \|\zeta_0 \Pi_j \zeta_0\|_2 \leq c_2 j^{-1/4}
     \ee
     with $c_2$ independent of $j$. Putting together \eqref{fin2},  \eqref{fin3}, and
     \eqref{fin4}, we conclude that there exists $c_3$ such that
     \bel{fin5}
     \left\|\sum_{j \in {\mathbb N}\, : \, B_j\in {\mathcal J}_+} \zeta_0 R_{\lambda}^-(B_j) {\mathcal W} \zeta_0 \Pi_j \zeta_0\right\|_1 \leq c_3 \sum_{j \in {\mathbb N}} j^{-5/4} < \infty.
     \ee
     Finally, we  estimate the trace-class norm of the second (finite) sum at the r.h.s. of \eqref{fin1}. We have
     $$
     \left\|\sum_{j \in {\mathbb N}\, : \, B_j\in {\mathcal J}_-} \zeta_0 R_{\lambda}^+(B_j) {\mathcal W}
     \zeta_0 \Pi_j \zeta_0\right\|_1 \leq \|\zeta_0\|_{L^{\infty}(\rd)}
     \sum_{j \in {\mathbb N}\, : \, B_j\in {\mathcal J}_-} \left\| R_{\lambda}^+(B_j) {\mathcal W}\right\| \| \zeta_0  \Pi_j \zeta_0\|_1.
     $$
     Moreover,
     \bel{fin10}
     \| \zeta_0  \Pi_j \zeta_0\|_1 = \|\Pi_j \zeta_0\|_2^2 = \frac{B}{2\pi} \|\zeta_0\|_{L^{2}(\rd)}^2, \quad j \in {\mathbb N},
     \ee
     (see e.g. \cite[Lemma 3.1]{FR}). Since the number of the Landau levels $B_j$ lying on ${\mathcal J}_-$ is finite,
     and the operators $R_{\lambda}^+(B_j) {\mathcal W}$ are bounded provided that $B_j \in {\mathcal J}_-$,  we get
     \bel{fin9}
     \left\|\sum_{j \in {\mathbb N}\, : \, B_j\in {\mathcal J}_-}
     \zeta_0 R_{\lambda}^+(B_j) {\mathcal W} \zeta_0 \Pi_j \zeta_0\right\|_1 < \infty.
     \ee
     Combining \eqref{fin5a} -- \eqref{fin7}, \eqref{fin8},
     \eqref{fin5}, and \eqref{fin9}, we find that the operator $P_{\lambda,{\mathcal E}} - P_{0,{\mathcal E}}$
     is trace-class.
     \end{proof}
\begin{lemma} \label{finl2}
    Let ${\bf a} \in C^1(\rd; \rd)$ with
    \bel{fin15}
    \|{\bf a}\|_{L^{\infty}(\rd)} + \|\nabla {\bf a}\|_{L^{\infty}(\rd)} \leq
    K,
    \ee
    with some $K < \infty$.
     Suppose that
    ${\mathcal E} \in \re$ is
    in a spectral gap of the family $H_\lambda$, $\lambda \in
    [0,\lambda_0]$ with some $\lambda_0 > 0$. Let $\phi_x$ be a
    smooth characteristic function of the unit square centered at
    $x \in \rd$. Then  for any $x, y
    \in \rd$ the operator $\phi_x P_{\lambda, {\mathcal E}}
    \phi_y$ is trace-class, and we have
    \bel{fin19}
    \|\phi_x P_{\lambda, {\mathcal E}}
    \phi_y\|_1 \leq C,
    \ee
    with $C = C(K)$ independent of $x, y \in \rd$, $\lambda \in
    [0,\lambda_0]$, and of ${\bf a}$ satisfying \eqref{fin15}.
\end{lemma}
The proof of the proposition is quite similar to the previous one
so that we omit the details, and only note that the analogues of
the bounds obtained in the proof of Lemma \ref{finl1} remain
uniform with respect to ${\bf a}$ satisfying \eqref{fin15}, and the
norms $\|\zeta_j\|_{L^p(\rd)}$ of the cut-off functions, as well
as the constant $c_2$ in \eqref{fin4}, are invariant under
translations of
the supports of $\zeta_j$.\\


{\large \bf Acknowledgements}. The authors thank Alex Sobolev for
providing the text of \cite{Be}.\\
N. Dombrowski and F. Germinet were supported in part by ANR 08
BLAN 0261. N. Dombrowski and G. Raikov were partially supported by {\em
N\'ucleo Cient\'ifico ICM} P07-027-F ``{\em Mathematical Theory of
Quantum and Classical Magnetic Systems"}. F. Germinet and G. Raikov were partially supported by
the Chilean Science Foundation {\em Fondecyt} under Grant 1090467.
G. Raikov thanks the Bernoulli Center, EPFL, Lausanne,
for a partial support during his participation in the Program
``{\em Spectral and Dynamical Properties of Quantum Hamiltonians}" in 2010.

\end{document}